\documentclass[11pt]{amsart}
\usepackage[margin=3cm]{geometry}               
\usepackage[usenames,dvipsnames,svgnames,table]{xcolor}    
\usepackage{graphicx}
\usepackage{tikz}
\usepackage{tikz-cd}
\usetikzlibrary{decorations.pathreplacing}
\usepackage{lineno, hyperref}
\usepackage{amssymb,amsfonts, amsmath, amsthm}
\usepackage{epstopdf}
\usepackage[shortlabels]{enumitem}
\usepackage{ulem}
\usepackage{cancel}
\usepackage{thmtools} 
\usepackage{mathtools} 
\usepackage{mathrsfs}
\usepackage{yhmath}

\title{Local quasi-isometries and tangent cones of definable germs}

\newtheorem{thm}{Theorem}[section]
\newtheorem{lem}[thm]{Lemma}
\newtheorem{prop}[thm]{Proposition}

\newtheorem{cor}[thm]{Corollary}
\newtheorem{example}[thm]{Example}
\newtheorem{rem}[thm]{Remark}

\numberwithin{equation}{section}

\newcommand{\bb}{\mathbb}
\newcommand{\al}{\mathcal}

\newcommand{\bff}{\mathbf}

\newcommand{\dist}{{\rm dist}}

\newcommand{\Seq}{{\rm Seq }}
\newcommand{\length}{{\rm length }}

\author[Nhan Nguyen]{Nhan Nguyen}
\address{FPT University, Danang, Vietnam}
\email{nguyenxuanvietnhan@gmail.com}

\begin{document}

\maketitle
\begin{abstract} In this paper, we introduce the notion of local quasi-isometry for metric germs and prove that two definable germs are quasi-isometric if and only if their tangent cones are bi-Lipschitz homeomorphic. Since bi-Lipschitz equivalence is a particular case of local quasi-isometric equivalence, we obtain Sampaio's tangent cone theorem as a corollary. As an application, we reprove the theorem by Fernandes-Sampaio, which states that the tangent cone of a Lipschitz normally embedded germ is also Lipschitz normally embedded.
\end{abstract}



\normalem 

\section{Introduction}
In this work, our focus is on the question of when the tangent cones of two given definable germs are bi-Lipschitz homeomorphic. Several authors have provided partial answers to this question. For instance, Bernig and Lytchak \cite{be-ly} demonstrated that if two subanalytic germs are sub-analytically bi-Lipschitz homeomorphic, then their tangent cones are also bi-Lipschitz homeomorphic. Sampaio \cite{Edson} showed that Bernig-Lytchak's result still holds if the germs are bi-Lipschitz equivalent, not necessarily subanalytic. Recently, Koike and Paunescu \cite{k-p} extended Sampaio's result by providing an analogous for arbitrary germs in $\mathbb{R}^n$ satisfying the sequence selection property (SSP). It is important to note that SSP is a necessary and sufficient condition for a germ in $\mathbb{R}^n$ to have a unique tangent cone (see \cite{S-S}).

Consider the semialgebraic sets $X = \{(x, y) \in \mathbb{R}^2 \mid y^2 = x^3, y \geq 0\}$ and $Y = \{(x, y) \in \mathbb{R}^2 \mid y^2 = x^5\}$, with $Z = X \cup Y$ (see Figure \ref{Fig1}). It can be easily verified that the bi-Lipschitz types (even the topological types) of the above sets at the origin are different. However, the tangent cones coincide and thus are bi-Lipschitz homeomorphic. This suggests that the bi-Lipschitz invariance of tangent cones should hold under weaker conditions on the given germs.

 \begin{figure}[h]
 \begin{center}
\begin{tikzpicture}

\draw[thick] plot [smooth,tension=1] coordinates{ (0,0) (1.2, 0.6 ) (2, 2) };
\draw (0,0) node[anchor=north] {$0$};
\draw (1,-2) node[anchor=south] {$X$};
\fill (0,0) circle (1.5pt);

\draw[thick, blue] plot [smooth,tension=1] coordinates{ (4,0) (5.3, 0.4 ) (6, 1.5) };
\draw[thick, blue] plot [smooth,tension=1] coordinates{ (4,0) (5.3, -0.4 ) (6, -1.5) };
\draw (4,0) node[anchor=north] {$0$};
\fill (4,0) circle (1.5pt);
\draw (5,-2) node[anchor=south] {$Y$};

\draw[thick] plot [smooth,tension=1] coordinates{ (8,0) (9.2, 0.6 ) (10, 2) };
\draw[thick, blue] plot [smooth,tension=1] coordinates{ (8,0) (9.3, 0.4 ) (10, 1.5) };
\draw[thick, blue] plot [smooth,tension=1] coordinates{ (8,0) (9.3, -0.4 ) (10, -1.5) };
\draw (8,0) node[anchor=north] {$0$};
\fill (8,0) circle (1.5pt);
\draw (9,-2) node[anchor=south] {$Z = X \cup Y$};

\end{tikzpicture}
\end{center}
\caption{}
  \label{Fig1}
\end{figure}

In the context of metric spaces, one has the notion of asymptotic tangent cones which are  defined as ultralimits of scaled pointed metric spaces $(X, d_X/n, p_i)$.  It is well-known that the bi-Lipschitz types of asymptotic tangent cones are quasi-isometry invariants of metric spaces. Quasi-isometry  plays a crucial role in geometric group theory following the fundamental work of Gromov \cite{Gromov}.  A map $f: (X, d_X)\to (Y, d_Y)$ between two metric spaces is called a {\bf quasi-isometry} if there are constants $K\geq 1$, $A\geq 0$ and $C\geq 0$ such that 

(i) $ \forall x, x'\in X, \frac{1}{K}d_X(x, x') - A \leq d_Y(f(x), f(x')) \leq K d_X(x, x') + A.$

(ii)  $\forall y \in Y, \exists x\in X: d_Y(f(x), y) \leq C$.

If $A = C = 0$ then $f$ is a { \bf bi-Lipschitz homeomorphism}. In particular, when $K = 1, A = C =0$ then $f$ is called an {\bf isometry}.

Note that quasi-isometry disregards the local geometry. To capture the local behavior, we propose the concept of \textbf{local quasi-isometry} (or quasi-isometry of metric germs). This notion is inspired by the definition of quasi-isometry. Since we work with germs, the constants $A$ and $C$ in conditions (i) and (ii) are replaced by germs of non-decreasing infinitesimal functions (see Section \ref{sec2}). We first prove that local quasi-isometry guarantees the invariance of the bi-Lipschitz type of blow-ups (see Theorem \ref{thm1}). Furthermore, if the given germs are subsets of Euclidean spaces and have unique tangent cones, then the blow-up of $X$ at a point $p$ is isometric to the usual tangent cone of $X$ at $p$. Consequently, if two such germs are quasi-isometric, then their tangent cones are bi-Lipschitz equivalent. We further establish that, in the o-minimal context, two definable germs are quasi-isometric if and only if their tangent cones are bi-Lipschitz homeomorphic (see Theorem \ref{main_theorem}). This provides a complete answer to the question mentioned at the beginning. In particular, we obtain Sampaio's tangent cone theorem as a corollary (Theorem \ref{thm_S}).

As an application, we recover the proof of the main theorem in \cite{Fer-Sam1}, which states that given a germ $(X, p)$ in $\mathbb{R}^n$ with a unique tangent cone, if  $(X, p)$ is Lipschitz normally embedded, then its tangent cone is also Lipschitz normally embedded (see Theorem \ref{thm_FS}).

The article is organized as follows:
\begin{itemize}
    \item Section \ref{sec2}: Introduces the concept of local quasi-isometry and explores its fundamental properties. 
    \item Section \ref{sec3}: Discusses the definition of blow-ups for metric germs and establishes the result that quasi-isometric metric germs yield bi-Lipschitz homeomorphic blow-ups.
    \item Sections \ref{sec4} and \ref{sec5}: Investigate tangent cones at specific points and at infinity for subsets of Euclidean spaces, respectively.
    \item Section \ref{sec6}: Examines the case of definable sets in o-minimal structures.
\end{itemize}

Throughout the note, we fix an o-minimal structure on $(\bb R,+, .)$. By “definable” we mean definable in the structure.  Given two non-negative functions $f, g: A \to \bb R$, we write $f\lesssim g$ if there is $C>0$ such that $f(x) \leq  C g(x)$ for every $x \in A$, and write $f\sim g$ if $f\lesssim g$ and $g\lesssim f$.

\section{Local quasi-isometries}\label{sec2}
\subsection{Quasi-isometry at a point}
{\bf A horn function at $0$} is the germ at $0$ of a non-decreasing function $\varphi: \bb R_{\geq 0} \to \bb R_{\geq 0} $ such that $$\lim_{t\to 0} \frac{\varphi(t)}{t} = 0.$$

Let $(X, p, d_X)$ and $(Y, q, d_Y)$ be two metric germs  and let $\varphi$ be a horn function at $0$.   A map $f: (X, d_X, p) \to  (Y, d_Y,  q)$ is called a {\bf $\varphi$-quasi-isometric embedding} if there are constants $L\geq 1$, $L_1, L_2>0$  and a neighborhood $U$ of $p$ in $X$ such that for all $x, x'$ in $U$ we have
\begin{equation}\label{equ_1}
    \frac{1}{L} d_X(x, x') - \delta(x,x', L_1, L_2) \leq d_Y(f(x), f(x')) \leq L d_X(x, x') + \delta(x,x', L_1, L_2),
\end{equation}
where $\delta(x,x', L_1, L_2) = L_1 \varphi \left(L_2 (d_X(x, p) + d_X(x', p))\right).$

We call $f$ a {\bf quasi-isometric embedding} if it is a $\varphi$-quasi-isometric embedding for some horn function $\varphi$.

The map $f$ is call a {\bf $\varphi$-quasi-isometry} if it is a $\varphi$-quasi-isometric embedding and there are a neighborhood $V$ of $q$ in $Y$ and constants constants $L_1', L_2'> 0$ such that for every $y\in V $ there is $x \in U$  such that 
\begin{equation}\label{equ_2}
    d_Y (y, f(x)) \leq L'_1 \varphi (L_2' d_X(x, p)).
\end{equation}

We say that  $(X, p)$ and $(Y, q)$ are  {\bf quasi-isometric} if there is a $\varphi$-quasi-isometry between them for some horn function $\varphi$.

A map $g: (Y, q) \to (X, p)$ is called a {\bf $\varphi$-quasi-inverse of $f$} if the following conditions are satisfied:

   ($\star$) there are a neighborhood $U'$ of $p$ in $X$ and constants $L_3, L_4 >0$ such that for all $x \in U'$
   $$ d_X(g\circ f (x),x ) \leq L_3 \varphi (L_4 d_Y (f(x), q))$$

   ($\star \star$) there are  a neighborhood $V'$ of $q$ in $Y$ and constants $L'_3, L'_4 >0$ such that for all $y \in V'$
   $$ d_Y(f\circ g(y), y) \leq L_3' \varphi (L'_4 d_X(g(y), p)).$$

   \begin{rem}\label{rem_1}\rm
   
    \begin{itemize}
    \item[(i)] A quasi-isometry needs neither to be continuous nor injective.
     \item[(ii)] $0$-quasi-isometry is the germ of a bi-Lipschitz homeomorphism.
    \item [(iii)] Let $f: (X, p) \to (Y, q)$ be a quasi-isometric embedding. Then,
            \begin{equation*} 
                 d_X(x, p) \sim d_Y(f(x), q) 
            \end{equation*}
   
    \end{itemize}
   \end{rem}

   \begin{example}\rm
   Let $X$, $Y$ and $Z$ be as in Figure \ref{Fig1}. Let $h: (X, 0) \to (Y, 0)$ be defined by $h(x, y) = (x, y^{5/3})$. It is easy to check that $h$ is a $\varphi$-quasi-isometry with $\varphi (t) = t^b$ where $1< b < \frac{3}{2}$, and the map $g: Y \to X$ defined by $g(x, y) = (x, |y|^{3/5})$ is a $\varphi$-quasi-inverse of $h$.  Therefore, $(X, 0)$ and $(Y, 0)$ are quasi-isometric. Similarly, one can show that $(X, 0)$ and $(Z, 0)$ are quasi-isometric as well. 
   \end{example}

\begin{lem}\label{lem_2} 

(i) Every $\varphi$-quasi-isometry has a $\varphi$-quasi-inverse.

(ii) Every quasi-inverse of a  $\varphi$-quasi-isometry is also a $\varphi$-quasi-isometry.

(iii) Composition of $\varphi$-quasi-isometries is again a $\varphi$-quasi-isometry.  
\end{lem}

\begin{proof} Let $f: (X, p)\to (Y, q)$ be a $\varphi$-quasi-isometry with some horn function (at $0$) $\varphi$.

(i) Let $g: (Y, q) \to (X, p)$ be a map defined as follows. Let $y \in Y$. If $y \in f(X)$, we take a point $x \in f^{-1} (y)$ and set $g(y) = x$. If $y\not \in f(X)$, by (\ref{equ_2}), there is $x\in X$ such that 
\begin{equation}\label{equ_3}
    d_Y(y, f(x)) \leq L_1'\varphi(L_2' d_X (x, p))
\end{equation}
where $L_1'$ and $L_2'$ are positive constants independent of $x$.
We then put $g(y) = x$.  We will show that $g$ is a $\varphi$-quasi-inverse of $f$. 

Let $x\in X$ near $p$. If $x \in g\circ f (X)$, then $x = g\circ f (x)$, so the condition ($\star$) is obvious. Assume that $x \not \in  g\circ f (X)$.  Set $x' = g\circ f(x)$ and $y = f(x)$. Note that $x'\neq x$ and $f(x) = f(x') = y$. It follows from (\ref{equ_1}) 
\begin{align*}
    \frac{1}{L} d_X(x, x') &\leq d_Y(f(x), f(x')) +   L_1 \varphi \left(L_2 \left(d_X(x, p) + d_X(x', p)\right)\right) \\
    &\leq  L_1 \varphi \left(L_2 \left(d_X(x, p) + d_X(x', p)\right)\right)
\end{align*}
for some constant $L, L_1, L_2$ independent of $(x, x')$.
By (iii) of Remark \ref{rem_1}, we can choose $L_3 >0$ (independent of $(x, x')$) such that $L_2 (d_X(x, p) + d_X(x', p)) \leq L_3 (d_Y(f(x), q) + d_Y(f(x'), q))$. 
Since $f(x) = f(x')$ and $\varphi$ is non-decreasing, 
$$ d_X(x, x')  \leq L L_1 \varphi (2L_3 d_Y(f(x), q)).$$
Thus, the condition ($\star$) is satisfied. 

We now check the condition ($\star \star$). Let $y \in Y$ near $q$. If $y \in f(X)$, by the construction, $f\circ g(y) = y$, hence ($\star \star$) is obvious. Assume that $y \not \in f(X)$. Let $x = g(y)$. It follows from (\ref{equ_3}), 
$$ d_Y(y, f\circ g(y)) = d_Y (y, f(x)) \leq L_1'\varphi(L_2' d_X(x, p)) = L_1'\varphi(L_2' d_X(g(y), p)).$$ This yields ($\star \star$). Therefore, $g$ is a $\varphi$-quasi-inverse of $f$.

(ii) Suppose $g: (Y, q) \to (X, p)$  is an $ \varphi$-quasi-inverse of $f$. We will show that $g$ is a $\varphi$-quasi-isometry. We need to show that $g$ satisfies (\ref{equ_1}) and (\ref{equ_2}).

Let $y, y'$ be points in $Y$ near $q$. Since $f$ is a $\varphi$-quasi-isometry,  by (\ref{equ_1}),
$$d_X(g(y),g( y')) \leq L d_Y(f \circ g(y), f\circ g(y')) + \delta(g(y), g(y'), L_1, L_2).$$
Moreover, 
\begin{align*}
    d_Y(f \circ g(y), f\circ g(y')) & \leq d_Y(y, y') +  d_Y(f \circ g(y), y) + d_Y(f \circ g(y'), y')  \\
    & \leq d_Y(y, y') +  C \varphi (C' d_Y(g(y), p)) + C \varphi (C' d_Y(g(y'), p)) \hspace{0.3 cm} \text{(by } (\star \star)) 
    \\
    &\leq d_Y(y, y') +  C \varphi\left(C'\left[d_Y(g(y), p) + d_Y(g(y'), p)\right]\right) \hspace{0.3 cm} \text{(since } \varphi \text{ is non-decreasing)}
\end{align*}
Here $C, C'$ are positive constants independent of $(y, y')$.
Recall that 
$$\delta(g(y), g(y'), L_1, L_2) = L_1 \varphi\left(L_2\left[d_Y(g(y), p) + d_Y(g(y'), p)\right]\right).$$
Hence, 
$$ d_X(g(y),g( y')) \leq L d_Y(y, y') + C_1 \varphi \left(C_1'\left[d_Y(g(y), p) + d_Y(g(y'), p)\right]\right)$$
where $C_1 = C + L_1$ and $C_1' = C' + L_2$.
This shows that $g$ satisfies the right-hand-side inequality in (\ref{equ_1}). By similar arguments, it is easy to check the left-hand-side inequality. 

Now let $z \in X$ near $p$. Take $y = f(z)$. Since $g$ is quasi-inverse of $f$, by ($\star$), there are constant $C, C'>0$ (independent of $y$) such that 
$$d_X(z, g (y))= d_X(z, g\circ f(z)) \leq C \varphi( C' d_Y(f(z), q)) =  C \varphi( C' d_Y(y, q)).$$
Thus, (\ref{equ_2}) is satisfied.

(iii) Suppose that $h: (Y, q) \to (Z, a)$ is a $\varphi$-quasi-isometry. We need to show that the map $h\circ f: (X, p) \to (Z, a)$ satisfies (\ref{equ_1}) and (\ref{equ_2}). We prove the right-hand-side inequality in (\ref{equ_1}), the left-hand-side one can be done similarly. Indeed, let $x, x' \in X$ near $p$. Since $f$ is a $\varphi$-quasi-isometric embedding, 
$$ d_X(x, x')   \leq L d_Y(f(x), f(x') + \delta(x, x', L_1, L_2) $$
where $L\geq 1, L_1, L_2 >0$ are constants independent of $(x, x')$.
In addition,  since $h$ is  a $\varphi$-quasi-isometric embedding, there are $C\geq 1$, $C_1, C_2 >0$ such that 
$$ d_Y(f(x), f(x') \leq C d_Z (h\circ f(x), h\circ f(x')) + \delta(f(x), f(x'), C_1, C_2).$$
This yields that 
$$ d_X(x, x')   \leq  LC  d_Z (h\circ f(x), h\circ f(x')) +  \delta(x, x', L_1, L_2)  + L \delta(f(x), f(x'), C_1, C_2).$$
By (iii) in Remark \ref{rem_1}, we can find $K>0$ big enough such that 
$$ K L_2(d_X (x, p) + d_X(x', p)) \geq C_2 (d_Y (f(x), q) + d_Y(f(x'), q)).$$
We may assume that $K L_1 >L C_1$. This implies that
$$ \delta(x, x', KL_1, KL_2)  \geq  L \delta(f(x), f(x'), C_1, C_2)$$ (since $\varphi$ is non-decreasing). Hence, 
$$ d_X(x, x')  \leq  LC  d_Z (h\circ f(x), h\circ f(x'))  + \delta (x, x', 2KL_1, KL_2),$$
so (\ref{equ_1}) is satisfied.

It remains to  check (\ref{equ_2}). Given $z \in Z$ near $a$, applying (\ref{equ_2}) to $h$, there is $y\in Y$ such that 
$$ d_Z (z, h(y)) \leq C \varphi(C'd_Y (y, q))$$
where $C, C'>0$ are constant independent of $z$. Let $g$ be a $\varphi$-quasi-inverse of $f$ and take $x = g(y)$. By (ii) in Remark \ref{rem_1},  we can choose $K>0$ depending only on $g$ such that $C'd_Y(y, q) \leq Kd_X(x, p)$. Since $\varphi$ is non-decreasing,  
$$ d_Z (z, h(y)) \leq C \varphi(Kd_X (x, p)).$$
Thus, $h\circ f$ satisfies (\ref{equ_2}). The proof is completed.    
\end{proof}

\section{Blow-ups}\label{sec3}
\subsection{Ultrafilters} Let $I$ be an infinite set. 
A {\bf filter} on $I$ is a non-empty family $\omega$ of subsets of $I$ such that 
\begin{enumerate}
    \item $\emptyset \not\in \omega$,
    \item If $A\in \omega$, $A\subset B$ then $B \in \omega$,
    \item If $A_1, \ldots, A_n$ are in $\omega$ then $A_1\cap \ldots \cap A_n \in \omega$.
\end{enumerate}
An element of $\omega$ is called {\bf $\omega$-large}. We say that a property $(P)$ holds for $\omega$-large if $(P)$ is satisfied for all in some $\omega$-large set.  A filter is called {\bf principal} if it contains a finite set.

An {\bf ultrafilter} is a filter $\omega$ satisfying further property that for any $A\in I$ either $A$ or $I \setminus A$ is in $\omega$. It is known that every filter contains an ultrafilter. In this paper, we only consider \textbf{non-principal ultrafilters}, i.e., ultrafilters that contain no finite sets. 

\subsection{Ultralimits}

Let   $\omega$ be a non-principal ultrafilter on $I$ and let $Y$ be a metric space. A {\bf $\omega$-limit} of a function 
$f: I \to Y$, denoted by 
$\lim_{\omega} f = \lim_{\omega} f(i)\in Y$, 
is defined to be a point $y \in Y$ if for every neighborhood $U$ of $y$ in $Y$, $f^{-1}(U)$ is $\omega$-large. 

An important basic fact is that if the $\omega$-limit exists then it is unique. In particular, if $Y$ is compact the $\omega$-limit is well-defined, meaning it uniquely exists. Consequently, if $f(i)$ is a bounded sequence in $\bb R^n$, then $\lim_\omega f(i)$ is uniquely determined.

To see this, let $y_1$ and $y_2$ be two different $\omega$-limits of $f$. Since $Y$ is a Hausdorff, there are neighbourhoods $U_{1}$ of $y_1$  and $U_{2}$ of $y_2$ in $Y$ such that $f^{-1}(U_{1}) \cap f^{-1}(U_{2}) = \emptyset$.  By definition, $f^{-1}(U_{1})$ and $f^{-1}(U_{2})$ both are in $\omega$, so is their intersection (see Property (3)). In other words, $\emptyset \in \omega$, which contradicts Property (1). Now suppose that $Y$ is compact. Assume on the contrary that the $\omega$-limit of $f$ does not exist, i.e., for any point $y\in Y$ there is a neighbourhood $U_y$ such that $f^{-1}(U_y) \not\in \omega$. By the compactness, $Y$ can be covered by finitely many such neighborhoods, let's say $U_{y_1}, \ldots, U_{y_k}$. By definition, $I \setminus f^{-1}(U_{y_i}) \in \omega $ for every $i = 1, \ldots, k$. This implies that $$\bigcap_i (I \setminus f^{-1}(U_{y_i})) = I\setminus \bigcup_i f^{-1}(U_{y_i}) \in \omega,$$
which is a contradiction since $I\setminus \bigcup_i f^{-1}(U_{y_i}) =  \emptyset$.

\subsection{Ultralimit of pointed metric spaces}
Let $\omega$ be a non-principal ultrafilter on $\bb N$. 
Let $(X_i, d_i, p_i)_{i \in \bb N}$ be a sequence of pointed metric spaces. Let $Seq$ denote the space of all sequences $(x_i)$ of points $x_i \in X_i$ such that $\sup\{ d_i(x_i, p_i)\} <\infty$. Given two sequences $\bff x = (x_i)$ and $\bff x' = (x'_i)$ in $Seq$, we define $d_\omega(\bff x, \bff x') = \lim_\omega d_i(x_i, x'_i)$. 
The ultralimit  $\lim_\omega(X_i, d_i, p_i)$ is the metric space $$(X_\omega, d_\omega) = (Seq, d_\omega)/ \sim$$ where $\bff x \sim \bff x'$ if $d_\omega (\bff x, \bff x') = 0$. We denote by $[(x_i)]$ the equivalence class of $(x_i)$.

\subsection{Blow-ups} Let $(X, d)$ be a metric space and let $p \in X$.  Let $\Lambda = (\lambda_i)_{i\in \bb N}$ be a sequence of positive real numbers  tending to infinity. The {\bf blow-up of $X$ at $p$} with the scale $\Lambda$ is the metric space 
$$(X_{\omega, \Lambda, p}, d_\omega) = \lim_\omega (X, d_i, p)$$ 
where $d_i = \lambda_i d$.

If $(X, d)$ is a unbounded metric space then we can define the {\bf blow-up of $X$ at infinity} to be the metric space 
$$(X_{\omega, \Lambda, \infty}, d_\omega) = \lim_\omega(X, d'_i, p)$$
where $d_i' = d/\lambda_i$.

\begin{thm}\label{thm1} Let $(X, p, d_X)$ and $(Y, q, d_Y)$  be two pointed metric spaces. If they are quasi-isometric then their blow-ups
$X_{\omega, \Lambda, p}$ and $Y_{\omega, \Lambda, q}$ are bi-Lipschitz homeomorphic. 
\end{thm}

\begin{proof} Let $f: (X, p) \to (Y, q)$ be a $\varphi$-quasi-isometry where $\varphi$ is a horn function at $0$. Consider the map  $f_{\omega, \Lambda}: (X_{\omega, \Lambda, p}, d_{X, \omega}) \to (Y_{\omega, \Lambda, q}, d_{Y, \omega})$ defined by 
$f([(x_i)]) = [(f(x_i))]$. We will show that $f_{\omega, \Lambda}$  is a bi-Lipschitz homeomorphism.

Indeed, let $\bff x = [(x_i)]$ and $\bff x' = [(x'_i)]$ be in  $X_{\omega, \Lambda, p}$,  by definition, 
$$ d_{Y,\omega}(f_{\omega, \Lambda}(\bff x), f_{\omega, \Lambda}(\bff y)) =  \lim_\omega \lambda_i d_Y (f(x_i), f(x'_i)).$$
It follows from (\ref{equ_1}) that there are constants $L\geq 1$, $L_1, L_2>0$ independent of $(\bff x, \bff x')$ such that 
$$ \frac{1}{L} d_X(x_i, x'_i) - \delta (x_i, x'_i, L_1, L_2) \leq d_Y(f(x_i), f(x'_i)) \leq L d_X(x_i, x'_i) + \delta (x_i, x'_i, L_1, L_2)$$
Hence, 
$$ \frac{1}{L}\lambda_i d_X(x_i, x'_i) - \lambda_i\delta (x_i, x'_i, L_1, L_2) \leq \lambda_i d_Y(f(x_i), f(x'_i)) \leq L \lambda_id_X(x_i, x'_i) + \lambda_i\delta (x_i, x'_i, L_1, L_2)$$
Observe that 
\begin{align*}
    \lambda_i\delta (x_i, x'_i, L_1, L_2) &= \lambda_i L_1 \varphi (L_2 [d_X (x_i, p) + d_X (x_i', p)])\\
& = L_1 L_2 \lambda_i [d_X (x_i, p) + d_X(x_i', p) ]\frac{\varphi (L_2 [d_X (x_i, p) + d_X (x_i', p)])}{L_2 [d_X (x_i, p) + d_X (x_i', p)]} \to 0 \text{ as } i \to \infty
\end{align*}
(since $\lambda_i (d_X (x_i, p) + d_X(x_i', p))$ is bounded  and $\varphi (x)/x \to 0$ as $x \to 0$). Passing the ultralimit, we get 
$$\frac{1}{L} d_{X,\omega} (\bff x, \bff x')\leq  d_{Y, \omega }(f_{\omega, \Lambda }(\bff x), f_{\omega, \Lambda}(\bff x')) \leq  L  d_{X,\omega} (\bff x, \bff x').$$
This shows that $f_{\omega, \Lambda }$ is a bi-Lipschitz embedding. 

To prove $f_{\omega, \Lambda }$ is a bi-Lipschitz homeomorphism, it suffices to show that $f_{\omega, \Lambda }$ is surjective. Let $\bff y = [(y_i)] \in Y_{\omega, \Lambda, q}$. By (\ref{equ_2}), there is a sequence $(x_i) \subset X$ such that 
\begin{equation}\label{equ_5}
    d_Y(y_i, f(x_i)) \leq C \varphi (C' d_X (x_i, p)).
\end{equation} 
for some constants $C, C'>0$ independent of $(x_i)$.

It follows that  
 \begin{equation}\label{equ_6}
     \frac{d_Y(y_i, f(x_i))}{ d_X(x_i, p)} \leq C\frac{\varphi (C' d_X (x_i, p))}{ d_X(x_i, p)} \to 0 \text{ as } i \to \infty.
 \end{equation}
Since  $d_Y(y_i, q) \geq   d_Y (f(x_i), q) - d_Y(f(x_i), y_i)$, by (\ref{equ_6}) and the fact that $d_Y (f(x_i), q) \sim d_X (x_i, p)$,
$$d_Y(y_i, q)\sim  d_Y(f(x_i), q) \sim d_X (x_i, p).$$
Moreover,  $\lambda_i d_Y(y_i, q)$ is bounded, so is $\lambda_i d_X(x_i, p)$, and hence $(x_i) \in \Seq_{\omega, \Lambda}(X,p)$. 

It follows again from (\ref{equ_5}) that 
 $$\lambda_i d_Y(y_i, f(x_i)) \leq C C' \lambda_i d_X(x_i, p) \frac{\varphi (C' d_X (x_i, p))}{ C' d_X(x_i, p)} \to 0 \text{ as } i \to \infty$$
where $C, C'$ are positive constants.
This implies that $d_{Y, \omega}(\bff y, f_{\omega, \Lambda } (\bff x))= 0$. In other words, $\bff y = f_{\omega, \Lambda } (\bff x)$) where $\bff x = [(x_i)]$. Thus, $f_{\omega, \Lambda }$ is surjective. 
\end{proof}

\section{Tangent cones at a point $p$}\label{sec4}
Let $X$ be a subset of $\bb R^n$ and let $p\in \overline{X}$. Fix a sequence $\Lambda = (\lambda_i)$ such that $\lim_{i \to \infty} \lambda_i = \infty$.
The \textbf{$\Lambda$-tangent cone} of $X$ at $p$ is the set 
$$ C (X, \Lambda, p) = \{ v \in \bb R^n: \exists (x_i) \in X, x_i \to p, \lim_{i \to \infty} \lambda_i (x_i-p)  = v\}.$$

When the set $C(X, p, \Lambda)$ is independent of the choice of $\Lambda$ then we say that $X$ has  \textbf{a unique tangent cone} at $p$.  

Set
$$C (X, p) =  \{ v\in \bb R^n, \exists (x_i) \subset X, x_i \to p, \exists (\lambda_i) \subset \bb R_+, \lambda_i \to \infty,  v = \lim_{i\to \infty}\lambda_i (x_i-p)\}$$
Then, if $X$ has a unique tangent cone at $p$ then 
$$ C(X, \Lambda, p) = C(X, p)$$ for all sequences $\Lambda$ tending to $\infty$.

Now let us assume that $X$ has a unique tangent cone at $p$. Fix a sequence  $\Lambda = (\lambda_i)$ as above and  fix  a principal ultrafilter $\omega$ on $\bb N$. 

Consider the following natural map: 
$$\Phi: (X_{\omega, \Lambda, p}, d_{out,\omega}) \to (C (X, p), d_{out}), \hspace{0.5cm} \Phi([(x_i)]) = \lim_\omega \lambda_i (x_i-p)$$
where $d_{out}$ is the metric induced by the Euclidean distance.

Since $(\lambda_i (x_i-p))_{i \in \bb N}$ is a bounded sequence in $\bb R^n$, the $\omega$-limit is uniquely determined. Moreover,  given $v = \lim_\omega \lambda_i (x_i-p)$, there is subsequence $(\lambda_{i_k}( x_{i_k}-p))$ such that $v$ is the limit of the sequence  $(\lambda_{i_k} (x_{i_k} -p))$ in the standard sense. This implies that $v \in C(X, p)$ and hence $\Phi$ is well defined. 

\begin{prop}\label{pro1}

$\Phi$ is an isometry.
\end{prop}
\begin{proof}
It suffices to show that $\Phi$ is a surjective isometric embedding. Without loss of generality, we may assume that $p$ coincides with the origin $0$. 
Let $\bff x = [(x_i)]$ and $\bff y = [(y_i)]$ be  in $X_{\omega, \Lambda, p}$. We have

$$d_{out,\omega}(\bff x, \bff y)   = \lim_\omega \lambda_i \| x_i - y_i\| = \|\lim_\omega \lambda_i x_i  - \lim_\omega \lambda_i y_i\| = \|\Phi (\bff x) - \Phi (\bff y)\|.$$
This implies that $\Phi$ is an isometric embedding. 

Let us now prove the surjectivity. Let $v\in C(X, p)$. By the uniqueness of the tangent cone, we have 
$$ C(X, p) = C(X, \Lambda, p)$$
It follows that there is a sequence $(x_i)$ on $X$ such that $\lim_{i \to \infty} \lambda_i x_i = v$. Put $\bff x = [(x_i)]$. It is obvious that $\bff x \in X_{\omega, \Lambda, p}$ and $\Phi(\bff x) = v$. Hence, $\Phi$ is surjective. 
\end{proof}
A direct consequence of Proposition \ref{pro1} and Theorem \ref{thm1} is as follow:

\begin{thm}\label{cor_main}
Let $(X, p)$ and $(Y, q)$ be germs in $\bb R^n$ and $\bb R^m$ respectively. Assume that $(X, p)$ and $(Y, q)$ have unique tangent cones. If $(X, p)$ and $(Y, q)$ are quasi-isometric then $C(X, p)$ and $C(Y, q)$ are bi-Lipschitz homeomorphic.
\end{thm}

\section{Quasi-isometries and tangent cones at infinity}\label{sec5}

\subsection{Quasi-isometry at infinity} {\bf A horn function at infinity} is the germ at infinity of a  non-decreasing function $\psi: \bb R_{\geq 0} \to \bb R_{\geq 0} $ satisfying $$\lim_{t\to \infty} \frac{\psi(t)}{t} = 0.$$ 

Let $(X, d_X)$ and $(Y, d_Y)$ be two unbounded metric spaces, and let $\psi$ be a horn function at infinity. A map $f: X \to Y$ is called a {\bf $\psi$-quasi-isometric embedding at infinity} if there are a point $p \in X$ and constants  $L\geq 1$, $L_1, L_2, r>0$ such that for all $x, x'$ in $X \setminus \bff B_r(X, p)$ we have
\begin{equation}\label{equ_1_1}
    \frac{1}{L} d_X(x, x') - \sigma(x,x', L_1, L_2) \leq d_Y(f(x), f(x')) \leq L d_X(x, x') + \sigma(x,x', L_1, L_2),
\end{equation}
where $\sigma(x,x', L_1, L_2) = L_1 \psi \left(L_2 (d_X(x, p) + d_X(x', p))\right)$ and $\bff B_r(X, p)$ is the ball in $X$ centered at $p$ of radius $r$.

We call $f$ a {\bf quasi-isometric embedding at infinity} if it is a $\psi$-quasi-isometric embedding at infinity for some horn function at infinity $\psi$.

The map $f$ is call a {\bf $\psi$-quasi-isometry at infinity} if it is a $\psi$-quasi-isometric embedding at infinity and there are  constants $L_1', L_2', r'>0$ such that for every $y\in Y\setminus \bff B_{r'} (Y, f(p))$ there is $x \in X$ such that  
\begin{equation}\label{equ_2_1}
    d_Y (y, f(x)) \leq L'_1 \psi (L_2' d_X(x, p)).
\end{equation}

We say that  $X$ and $Y$  are  {\bf quasi-isometric at infinity} (or $(X, \infty)$ and $(Y, \infty)$ are quasi-isometric) if there is a $\psi$-quasi-isometry between them for some horn function (at infinity) $\psi$.

Given $p$ as above, we set $q = f(p)$. A map $g: Y \to X $ is called {\bf $\psi$-quasi-inverse of $f$} if the following conditions are satisfied:

   ($\star'$) there are a neighborhood $U'$ of $\infty$ in $X$ and constants $L_3, L_4 >0$ such that for all $x \in U'$
   $$ d_X(g\circ f (x),x ) \leq L_3 \psi (L_4 d_Y (f(x), q))$$

   ($\star' \star'$) there are  a neighborhood $V'$ of $\infty$ in $Y$ and constants $L'_3, L'_4 >0$ such that for all $y \in V'$
   $$ d_Y(f\circ g(y), y) \leq L_3' \varphi (L'_4 d_X(g(y), p)).$$

   Here by a neighborhood at $\infty$ of a metric space $X$ (res. $Y$) we mean the set $X \setminus \bff B_R (X, p)$ (resp. $X \setminus \bff B_{R'} (Y, q)$) where $R$ and $R'$ are big enough real positive numbers.

\begin{rem}\rm

The definition of quasi-isometry at infinity is independent of the choice of $p$ that is if the map $f$ is a quasi-isometry, then for every $p\in X$, the inequalities (\ref{equ_1_1}) and (\ref{equ_2_1}) hold (the constants $L, L_1, L_2, L_1', L_2'$ may vary depending on the chosen $p$).

\end{rem}
By employing a similar reasoning as in the proof of Lemma \ref{lem_2}, we obtain 
\begin{lem}\label{lem_4_2} 
(i) Every $\psi$-quasi-isometry at infinity has a $\psi$-quasi-inverse.

(ii) Every quasi-inverse of a  $\psi$-quasi-isometry at infinity is also a $\psi$-quasi-isometry at infinity.

(iii) Composition of $\psi$-quasi-isometries at infinity is again a $\psi$-quasi-isometry at infinity.  
\end{lem}

\begin{thm}\label{thm2}  Let $(X, d_X)$ and $(Y, d_Y)$  be two unbounded  metric spaces. Then, if $(X, \infty)$ and $(Y, \infty)$ are quasi-isometric then their blow-ups at infinity
$(X_{\omega, \Lambda, \infty}, d_{X, \omega})$ and $(Y_{\omega, \Lambda, \infty}, d_{Y, \omega})$ are bi-Lipschitz homeomorphic. 
\end{thm}
\begin{proof}
    Let $f: (X, \infty) \to (Y, \infty)$  be a $\psi$-quasi-isometry where $\psi$ is a horn function at infinity. Consider the map 
    $$f_{\omega, \Lambda} : X_{\omega, \Lambda, \infty} \to Y_{\omega, \Lambda, \infty}$$
    defined by $[(x_i)]\mapsto [(f(x_i))]$. By the same arguments as in the proof of Theorem \ref{thm1}, it is easy to check that $f_{\omega, \Lambda}$ is a bi-Lipschitz homeomorphism.
\end{proof}

\subsection{Tangent cones at infinity}
Let $X$ be an unbounded subset of $\bb R^n$. Let $\Lambda = (\lambda_i) \subset \bb R$ be a sequence such that $\lim_{i \to \infty} \lambda_i = \infty.$ Similar to the definition of tangent cones at a point, the \textbf{$\Lambda$-tangent cone of $X$ at infinity} is the set 
$$ C (X, \Lambda, \infty) = \{ v \in \bb R^n: \exists (x_i) \in X, x_i \to \infty, \lim_{i \to \infty} \frac{x_i}{\lambda_i}  = v\}.$$

As the set $C(X,  \Lambda, \infty)$ is independent of the choice of $\Lambda$ then we say that $X$ has  \textbf{a unique tangent cone at infinity}.  

Set
$$C (X, \infty) =  \{ v\in \bb R^n, \exists (x_i) \subset X, x_i \to p, \exists (\lambda_i) \subset \bb R_+, \lambda_i \to \infty,  v = \lim_{i\to \infty}\frac{x_i}{\lambda_i}\}.$$
Then, if $X$ is of a unique tangent cone at infinity then 
$$ C(X, \Lambda, \infty) = C(X, \infty)$$ for every sequence $\Lambda$ tending to $\infty$.

\begin{prop}\label{pro2} Let $X$ be an unbounded subset of $\bb R^n$. Consider $X$ as a metric space with the outer metric. If $X$ is of a unique tangent cone at infinity, then the  map
     $$\Psi:  (X_{\omega, \Lambda, \infty}, d_{out,\omega}) \to (C (X, \infty), d_{out}), \hspace{0.5 cm} [(x_i)]\mapsto \lim_\omega \frac{x_i}{\lambda_i}$$
    is an isometry.
\end{prop}
\begin{proof}
    Similar to the proof of Proposition \ref{pro1}.
\end{proof}

Consequently, we get

\begin{cor}
    Let $X\subset \bb R^n$ and $Y\subset \bb R^m$ be unbounded sets of unique tangent cones at infinity. If $(X, \infty)$ and $(Y,\infty)$ are quasi-isometric then $C(X, \infty)$ and $C(Y, \infty)$ are bi-Lipschitz homeomorphic.
\end{cor}

\section{Applications}\label{sec6}

\subsection{Tangent cones of definable germs}
In this section, we examine germs of sets definable in an o-minimal structure. It is known as a consequence of Curve Selection that such germs possess unique tangent cones. Our aim is to prove the following result:

\begin{thm}\label{main_theorem} Let $(X, p)\subset \bb R^n$ and $(Y, q) \subset \bb R^m$ be two definable germs. Then, $(X, p)$ and $(Y, q)$ are quasi-isometric  if and only if their tangent cones are bi-Lipschitz homeomorphism.  
\end{thm}

We need some preparation before giving the proof. 

Given $A \subset \bb R^n$ and a non-negative function $\varphi: \bb R_\geq 0\to \bb R_{\geq 0}$, we define 
$$ \al N_\varphi (A) = \{ x \in \bb R^n: \dist(x, A) \leq \varphi(\|x\|)\}.$$
\begin{lem}\label{lem_neiborhood2}
    Let $X$ be a definable set in $\bb R^n$ and $0\in \overline{X}$. Then, the germs $(X, 0)$ and $(C(X, 0), 0)$ are quasi-isometric.
\end{lem}
\begin{proof}
It follows from \cite[Proposition 4.15]{Loi1} there is a definable germ $\theta: (\bb R, 0) \to (\bb R, 0)$ which is odd, strictly increasing such that 
$$ (X, 0) \subset \al N_{\varphi} (C(X, 0))$$ 
where $\varphi(t) = t\theta(t)$. It is clear that $\varphi$ is a horn function at $0$. 

Consider the set $Z=  \{(x, y) \in X \times C(X, 0): y = \dist (x, C(X, 0)\}$.
Since $C(X, 0)$ is a closed set, $Z_x = \{y\in Z: (x, y) \in Z\}$ is nonempty for every $x\in X$. By Definable Choice (see \cite{Dries}, Ch 6, (1.2)), there is a definable map $f: X \to C(X, 0)$ such that the graph of $f$ is contained in $Z$.

We will show that $f$ is the desired isometry.  It suffices to check that  $f$ satisfies the conditions (\ref{equ_1}) and (\ref{equ_2}). 

First note that 
$$\|f(x) - x\| \leq \varphi(\|x\|)$$ 
for every $x$ near the origin. Given $x, y$ near $0$, we have 
\begin{align*}
    \|f(x) - f(y)\| &\leq \|x - y\| + \|f(x) -x\| + \|f(y) - y\|\\
    & \leq \|x - y\| + \varphi (\|x\|) + \varphi(\|y\|)\\
    & \leq \|x - y\| + \varphi(\|x\| +\|y\|) \hspace{0.3cm} \text{(since } \varphi \text{ is increasing)}.
\end{align*}
Similarly, 
\begin{align*}
    \|x - y\| &\leq \|f(x) - f(y)\| + \|f(x) -x\| + \|f(y) - y\|\\
    & \leq \|f(x) - f(y)\| + \varphi (\|x\|) + \varphi(\|y\|)\\
    & \leq \|f(x) - f(y)\| + \varphi(\|x\| +\|y\|) \hspace{0.3cm} \text{(since } \varphi \text{ is increasing)}.
\end{align*}
The inequality (\ref{equ_1}) is then satisfied.

Now given $y \in (C(X, 0), 0)$, since $X \subset \al N_{\varphi} (C(X, 0))$, there is $x \in (X, 0)$ such that 
$\|x - y\| \leq \varphi(\|x\|)$. Then, 
$$ \|y - f(x)\| \leq \| y - x\| + \| x - f(x)\| \leq 2 \varphi(\|x\|).$$
Thus, the condition (\ref{equ_2}) is satisfied and the proof is complete.
\end{proof}

Now we are in a position to prove Theorem \ref{main_theorem}. 

{\bf Proof of Theorem \ref{main_theorem}}. The necessity follows from Theorem \ref{cor_main}. It remains to check the sufficiency. Without loss of generality, we may assume that $p$ and $q$ coincide with the origin $0$. From the assumption, $(C(X, 0), 0)$ and $(C(Y, 0), 0)$ are quasi-isometric. By Lemma \ref{lem_neiborhood2} we know that $(X, 0)$ and $(C(X, 0), 0)$ (resp. $(Y, 0)$ and $(C(Y, 0), 0)$) are quasi-isometric. Therefore, $(X, 0)$ and $(Y, 0)$ must be quasi-isometric (see Lemma \ref{lem_2}, property (iii)). $\qed$

\vspace{0.5 cm}

A direct consequence of Theorem \ref{main_theorem} is as follows: 
\begin{thm}[\cite{Edson}, Theorem 2.2]\label{thm_S}
    Let $X$ and $Y$ be two definable sets and let $p \in \overline{X}$ and $q \in \overline{Y}$. If $(X, p)$ and $(Y, q)$ are bi-Lipschitz homeomorphic, then  $C(X, p)$ and $C(Y, q)$ are also bi-Lipschitz homeomorphic.
\end{thm}

A version at infinity of Theorem \ref{main_theorem}  is as follows. 
\begin{thm}\label{thm_4}
    Let $X\subset \bb R^n$ and $Y\subset \bb R^m$ be two unbounded definable sets. If  $(X, \infty)$ and $(Y, \infty)$ are quasi-isometric then $C(X, \infty)$ and $C(Y, \infty)$ are bi-Lipschitz homeomorphic.
\end{thm}

\begin{proof}
The proof follows the same arguments presented in the proof of Theorem \ref{main_theorem}, using Lemma \ref{lem_neighborhood_2} and Lemma \ref{lem_4_2} instead of Lemma \ref{lem_2} and Proposition 4.15 in \cite{Loi1}.

\begin{lem}\label{lem_neighborhood_2}
    Let $X\subset \bb R^n$ be an unbounded definable set. Then, there is a continuous horn function (at infinity) $\psi: \bb R_{\geq 0} \to \bb R_{\geq 0}$ such that
    $$ (X, \infty) \subset \al N_{\psi} (C(X, \infty)).$$
\end{lem}
\begin{proof}

  For $x \in X$ we set $h(x) = \frac{\dist(x, C(X, \infty))}{\|x\|}$. 

Claim: $h(x) \to 0$ when $x$ tends to $\infty$. 

If the claim fails, by Curve Selection, there are a constant $c>0$ and a continuous definable curve $\gamma: (a, +\infty)\to X$, $a>0$, $\lim_{t\to\infty} (\gamma(t)) = \infty$ such that 
$h(\gamma(t)) > c$ for all $t$. By Monotonicity (\cite{Dries} Ch 3 (1.2)), $\lim_{t \to \infty} \frac{\gamma(t)}{\|\gamma(t)\|} = v  \in C(X, \infty)$. This yields that $ u(t) = \|\gamma(t)\|v $ is also in $C(X, \infty)$. Then, 
$$h(\gamma(t)) = \frac{\dist (\gamma (t), C(X, \infty))}{\|\gamma(t)\|} \leq\frac{\dist (\gamma (t), u(t))}{\|\gamma(t)\|} \to 0$$
which is a contradiction. Therefore, the claim is true. 

Now consider the diffeomorphism: $$\sigma: \bb R^n \setminus\{0\} \to \bb R^n \setminus\{0\}, \sigma(x) = \frac{x}{\|x\|^2}.$$
Set $Z = \sigma (X\setminus\{0\})$ and $g: Z \to \bb R$, $g = h\circ \sigma^{-1}$. Since  $\lim_{z\to 0} g(z) = 0$, it is possible to extend $g$ to a continuous function $\tilde{g}$ defined on $\overline{Z}$ and $\tilde{g}(0) = 0$. Fix an $r > 0$. By \L ojasiewicz inequality (see \cite{Loi2, D-M}), there is an odd, strictly increasing continuous definable function $\theta: \bb R \to \bb R$, $\theta(0) = 0$ such that $\theta(\tilde{g} (z)) \leq \|z\|$ for all $z \in \overline{Z} \cap \bff B^n_r$. Then, $\tilde{g}(z) \leq \theta^{-1}(\|z\|)$.

Let $\psi_1: \bb R_{\geq 0}\to \bb R$ be a function defined by $\psi_1(t) = \theta (1/t)$. Note that $\lim_{t\to \infty} \psi_1(t) = 0$. Given $x \in X \setminus \bff B^n_{1/r}$, by setting $z = \sigma(x)$ we have 
$$h(x) = \tilde{g} (z) \leq \theta^{-1}(\|z\|) = \psi_1(\|x\|) \hspace{0.3cm} \text{(since } \|z\| = 1/\|x\|).$$
This implies on $X \cap \bff B^n_{1/r}$, $\dist(x, C(X, \infty)) \leq \psi_1(\|x\|)\|x\|.$ Put $\psi_2 (t) = t \psi_1 (t)$. Then $X \cap \bff B^n_{1/r}$ is contained in $\al N_{\psi_2}(C(X, \infty))$. It suffices to show that there is a horn function at infinity $\psi \geq \psi_2$ when $t$ sufficiently large. Indeed if $\psi_2 (t) \leq \ln (t)$ as $t$ big enough,  we take $\psi(t) = \ln(t)$ which is obviously a horn function at infinity.  Otherwise, we choose $\psi = \psi_2$. Note that $\ln (t) \leq  \psi_2(t) \leq t$ when $t$ is big enough. Since  $\ln(t)$ and $t$ are increasing functions, so is $\psi_2$ as $t$ is big enough. This implies that $\psi = \psi_2$ is a horn function at infinity. This completes the proof. 
\end{proof}
\end{proof}

\subsection{Tangent cones of a Lipschitz normally embedded germ}
Let $X$ be a path-connected subset of $\bb R^n$. Consider following two metrics naturally equipped on $X$: {\it the outer metric ($d_{out}$)} induced by the Euclidean distance of the ambient space $\bb R^n$ and {\it the inner metric ($d_{inn}$)} where the distance between two points in $X$ is defined as the infimum of the lengths of rectifiable curve in $X$ connecting these points.  $X$ is called {\it Lipschitz normally embedded} (or {\it LNE} for brevity) if these two metrics are equivalent, i.e., there is a constant $C\geq 1$ such that 
$$\forall x, y \in X, d_{inn} (x,y) \leq C d_{out} (x,y).$$
The infimum of such constants $C$ is referred to as the {\it LNE constant} of $X$. 

It is obvious that $X$ is LNE if and only if the identity map  $$ id: (X, d_{inn}) \to (X, d_{out})$$ is a bi-Lipschitz homeomorphism.

We say that $X$ is {\it LNE at a point $x_0\in \overline{X}$} (or {\it the germ $(X, x_0)$ is LNE}) if there is a neighborhood $U$ of $x_0$ in $\bb R^n$ such that $X\cap U$ is LNE.

\begin{thm}[ \cite{Fer-Sam1}, Theorem 2.6]\label{thm_FS}
    Let X be a path-connected subset of $\bb R^m$ and let  $x_0 \in \overline{X}$. Suppose that X has a unique tangent cone at $x_0$. If $X$ is Lipschitz normally embedded, so is $C(X, x_0)$.
\end{thm}
\begin{proof}
Without loss of generality, we may assume that $x_0$ coincides with the origin $0$. 
    Consider the identity map 
    $id: (X, 0, d_{inn}) \to (X,0, d_{out})$. Since $(X, 0)$ is LNE so the map $id$ is a bi-Lipschitz homeomorphism.
Fix a non-principal ultrafilter $\omega$  on $\bb N$ and a sequence of positive numbers $\Lambda = (\lambda_i)_{i \in \bb N}$ such that $\lim_{i\to \infty} \lambda_i = \infty$.  

By Theorem \ref{thm1} and Theorem \ref{thm_4}, we have the following diagram 
\begin{center}
\begin{tikzcd}
 (X_{\omega, \Lambda, 0}, d_{inn, \omega})  \arrow[rd] \arrow[r, "id_{\omega, \Lambda} "] &(X_{\omega, \Lambda, 0}, d_{out, \omega}) \arrow[d, "\Phi"]  \\
& (C(X,0), d_{out})
\end{tikzcd}
\end{center}
where $id_{\omega, \Lambda}$ is a bi-Lipschitz homeomorphism and $\Phi$ is a isometry. 
Set $\Psi: = \Phi \circ id_{\omega, \Lambda}$. Then $\psi$ is a bi-Lipschitz homeomorphism. 
Let $x, y$ be two points in $C(X, 0)$.  By definition, there are sequences 
$(x_i)$ and $(y_i)$ in $X$ tending to $0$ such that $\lim_{i \to \infty} \lambda_i y_i = x$ and $\lim_{i \to \infty} \lambda_i y_i = y$.
For each $i$, let $\gamma_i: [0,1] \to X$ be a shortest arc in $X$ connecting $x_i$ and $y_i$. Reparametrizing $\gamma_i$, we may assume that $\length (\gamma_i ([0,t]))  = t. \length(\gamma_i)$.

It is obvious from definition that $\Psi ([(x_i)]) = x$ and $\Psi ([(y_i)]) = y$. Moreover,  
$$\lim_{i\to \infty} \lambda_i \|x_i - y_i\| = \|x -y\|.$$
Thus, there is $N\in \bb N$ big enough such that for every $i \geq N$, 
$$  \lambda_i \|x_i - y_i\| \leq 2 \|x -y\|.$$
Since $X$ is LNE, 
\begin{equation}\label{equ7}
   \lambda_i \length(\gamma_i) \leq  2 C\|x-y\|, 
\end{equation} 
where $C$ is the LNE constant of $(X, 0)$.   

Consider the following map 
$$ \gamma : [0, 1] \to (X_{\omega, \Lambda, 0}, d_{inn, \omega}) \hspace{0.5cm} t \mapsto [(\gamma_i(t))].$$
For $t, t' \in [0,1]$, $t >t'$, we have
\begin{equation}
d_{inn, \omega} ([(\gamma_i(t))], [(\gamma_i(t'))])  = \lim_{\omega} (\lambda_i d_{inn} (\gamma_i(t), \gamma_i(t'))\\
     \leq  \lim_{\omega} \lambda_i  (t- t') \length(\gamma_i)    
\end{equation}
Combining with (\ref{equ7}), we get 
    \begin{equation*}
       \frac{}{} d_{inn, \omega} ([(\gamma_i(t))], [(\gamma_i(t'))]) \leq  2 C (t - t') \|x - y\|
    \end{equation*}
This implies that 
\begin{itemize}
    \item [(i)] $\gamma$ is a continuous map, hence it is a curve in $(X_{\omega, \Lambda, 0}, d_{inn, \omega})$ connecting $[(x_i)]$ and $[(y_i)]$.
    \item [(ii)] $\length (\gamma) \leq 2 C \| x - y\|$.
\end{itemize}

Now, let $\tilde{\gamma}$ be the image of $\gamma$ by the map $\Psi$. Then,  $\tilde{\gamma}$ is a curve in $X$ connecting $x$ and $y$. Moreover,
$$ d_{inn} (x, y) \leq \length(\tilde{\gamma}) \leq   K \length (\gamma) \stackrel{(\rm ii)}{\leq} 2K C  \|x - y\|$$
where $K$ is the bi-Lipschitz constant of $\Psi$.
This ends the proof.

\end{proof}

\bibliographystyle{siam}
\bibliography{Biblio}

\end{document}